\documentclass[11pt]{amsart}
\usepackage{geometry}                % See geometry.pdf to learn the layout options. There are lots.
\usepackage[parfill]{parskip}    % Activate to begin paragraphs with an empty line rather than an indent
\usepackage{amsmath}
\usepackage{graphicx}
\usepackage{amssymb}
\usepackage{epstopdf}
\newtheorem{theorem}{Theorem}[section]
\newtheorem{lemma}[theorem]{Lemma}

\title{Pfister's theorem fails in the free case}
\author{Martin Harrison}

\begin{document}

\begin{abstract}
Artin solved Hilbert's $17^{th}$ problem by showing that every positive semidefinite polynomial can be realized as a sum of squares of rational functions. Pfister gave a bound on the number of squares of rational functions:  if $p$ is a positive semi-definite polynomial in $n$ variables, then there is a polynomial $q$ so that $q^2p$ is a sum of at most $2^n$ squares.

As shown by D'Angelo and Lebl, the analog of Pfister's theorem fails in the case of Hermitian polynomials.  Specifically, it was shown that the rank of any multiple of the polynomial $\|z\|^{2d} \equiv (\sum_j |z_j|^2)^d$ is bounded below by a quantity depending on $d$.  Here we prove that a similar result holds in a free $\ast$-algebra.
\end{abstract}
\maketitle
\section{Introduction}
The aim of this section is to define the main objects and to review some related work.  We work in the real free $\ast$-algebra $\mathbb{R}\langle X, X^* \rangle$ generated by the $n$ noncommuting (NC) variables $X_1, \ldots, X_n$ and their \emph{adjoints} $X_j^*$.  After taking a representation we can think of these variables as real square matrices, and the $\ast$ function on $\mathbb{R}\langle X, X^* \rangle$ as the transpose operation.  In particular, $\ast$ respects addition and multiplication by scalars and is defined on monomials by $(X_{j_1} \cdots X_{j_k}) ^*= X_{j_k}^* \cdots X_{j_1}^*$ and $(X_j^*)^*=X_j$.  We use multi-indices $\alpha$, tuples of non-negative integers from $0$ to $2n$, to index monomials: $X^\alpha \equiv X_{\alpha_1}X_{\alpha_2} \cdots X_{\alpha_k}$.  $X^\emptyset$ is simply the empty word, denoted by $1$.  For $0<j\leq n$, we define $X_{j+n} \equiv X_j^*$.  We define conjugation and concatenation of multi-indices $\alpha$ and $\beta$ by the equations $X^{\alpha^*}=(X^\alpha)^*$ and $X^{\alpha \circ \beta}  = X^{\alpha} X^{\beta}$.

Evaluation of $p \in \mathbb{R}\langle X, X^* \rangle$ at a tuple $(M_1, \ldots, M_n)$ of square matrices of the same size is defined by the substitution of $M_j$ for $X_j$ and $M_j^T$ for $X_j^*$.

We say that $p \in \mathbb{R}\langle X, X^* \rangle$ is symmetric when $p^*=p$.  Such a polynomial $p$ is said to be \emph{matrix positive} if the matrix $p(M)$ is positive semidefinite (or \emph{PSD}) for every tuple $M$ of square matrices.  It was shown by Helton in \cite{HE} that every matrix positive polynomial is a sum of squares (\emph{SOS}).  The minimal number of squares required to express a matrix positive polynomial as a sum of squares is not known in general, although upper bounds are easy to obtain.  The question is open in the commutative case as well, and in both cases amounts to a problem of rank minimization. A great many types of rank minimization problems have been successfully attacked in recent years with \emph{semidefinite programming} techniques (see \cite{FAZ} for examples).  A complete characterization of conditions for success of the nuclear norm approach, or ``trace-heuristic'', in this context is not know, though Recht provided in \cite{DBLP:conf/cdc/RechtXH08} a probabilistic characterization of success for particular classes of rank minimization problems.

Optimization in certain quantum physics problems is done over feasible regions of operators on Hilbert spaces, and so NC variables are useful there. Several examples and a general framework for such problems are presented in \cite{PIR}, where the semidefinite programming relaxations of Lasserre are extended to the NC setting. Motivation for the study of NC polynomials from control theory is discussed in \cite{HE2}.

\section{polynomials, associated matrices and sums of squares}

To any symmetric polynomial $p \in \mathbb{R}\langle X, X^* \rangle$ we can associate a real, symmetric matrix $M$ with the property

\[ V^*MV=p\]

where $V^*=(X^{\alpha^*})_{|\alpha|\leq d}$, and $V$ is the column vector $(X^{\alpha})_{|\alpha|\leq d}$ (with the monomials in graded lexicographical order). The matrix $M$ is not unique, in fact the set of all such matrices (for a fixed p) forms an affine space which we will denote $\mathcal M_p$.

By the rank of $p$, we mean the minimum of $\text{rank}(M)$ over all $M\in \mathcal {M}_p$.  For a positive polynomial, this minimum is to be taken over only the PSD matrices.  The following lemma helps us obtain a lower bound on rank

\begin{lemma}
If $A$ is a symmetric matrix satisfying $V^*AV=0$, then the $(2n)^d\times (2n)^d$ lower right submatrix of $A$ is the zero matrix.
\end{lemma}

\proof Let $B$ denote the block in question, and $\hat V$ the tautological vector of just the monomials of degree $d$. Then $V^*AV=0$ implies that $\hat V^* B \hat V=0$ as well since the product $\hat V^* B \hat V$ yields exactly the degree $2d$ terms of the polynomial $V^*AV$. But the entries of $B$ are exactly the coefficients of the distinct monomials in $\hat V^*B\hat V$, hence $B$ is the zero matrix.

The lemma above shows that there is no freedom in choosing the block corresponding to the degree $2d$ terms of the polynomial. Since the rank of this block gives a lower bound on the rank of the whole matrix, taking the block to be the the $(2n)^d\times (2n)^d$ identity yields a polynomial with rank at least $(2n)^d$.
\subsection{Positive polynomials and sums of squares}
In the commutative case it is well-known that the cone of positive polynomials properly contains the SOS cone. Motzkin's polynomial $M(x,y)=1+x^2y^4+y^2x^4-3x^2y^2$ is the first known example of a positive polynomial outside the SOS cone, and was discovered decades after Hilbert proved the existence of such polynomials.

In contrast, the NC setting offers the nice result, proved by Helton in \cite{HE}, that any positive polynomial is a sum of squares.  Here, a square takes the form $f^*f$, so that obviously a sum of squares is positive in the sense defined above.  In order to understand the SOS representation of a positive polynomial, we use the matrix representation introduced above.  The following lemma leads us to the semidefinite programming formulation of the rank minimization problem.

\begin{lemma}
A polynomial $p$ is matrix positive exactly when it can be expressed $p=V^*MV$, with $M$ a PSD matrix.  The rank of p is exactly the minimum number of squares over all SOS representations of $p$.
\end{lemma}

The proof is straightforward.  It follows that the minimum number of squares for a positive $p$ is 

\[\min \hspace{.5cm}  \text{rank  }X\]
\[\hspace{1cm} \text{s.t.   } V^*XV=p,\]
\[\hspace{2cm} X \succeq 0\] 

which can be calculated efficiently (but not always accurately), by minimizing instead the trace of $X$.

As a simple example of this problem consider the polynomial $P=1+X^*X+XX^*$, clearly a SOS. The polynomial $P$ is a sum of $3$ squares, but can be expressed as a sum of $2$ squares (and no fewer).  To see why we parameterize the affine space $\mathcal M_P$ by the single parameter $t \in \mathbb{R}$.  As usual  $V=(1,X,X^*)^T$, and so $P=V^*V=X^*X+XX^*+1$. Defining 

\[M=\begin{pmatrix} 
 0 & 1 & -1\\
 1 & 0 &0\\
  -1 & 0& 0
\end{pmatrix}\]

we get $\mathcal M_P= \{I+tM|t\in \mathbb{R}\}$, and find the minimal SOS representation $$P=\Big(X+\frac{\sqrt{2}}{2}\Big)^*\Big(X+\frac{\sqrt{2}}{2}\Big) +\Big(X^*-\frac{\sqrt{2}}{2}\Big)^*\Big(X^*-\frac{\sqrt{2}}{2}\Big) $$
on the boundary of the region where $I+tM \succ 0$. Note that in this example trace is constant on $\mathcal M_P \cap PSD$, and that the given solution is obtained by maximizing $t$ over $\{t|I+tM \succeq 0\}$.

\section{The Examples}
Pfister's Theorem gives a bound on the number of rational functions in the SOS representation of a PSD polynomial. The bound is remarkable because it does not depend on the degree of the polynomial in question. D'Angelo and Lebl proved in \cite{DAN} that this result fails for Hermitian polynomials.  We'll show that it fails for noncommutative polynomials. The first theorem below is needed for the second. It is easy to check that the polynomial $S$ below has rank $(2n)^d$, but more is true. 
\begin{theorem}
Suppose that $q \in \mathbb{R}\langle X, X^*\rangle$ and define $S=\sum_{|\alpha|=d} X^{\alpha^*}X^\alpha$. Then $p=q^*Sq$ has rank at least $(2n)^d$.  Here, $(2n)^d$ is the dimension of $span\{X^\alpha\}_{|\alpha|=d}$
\end{theorem}

\proof Since $p$ is matrix positive, it is a sum of squares, and so we may write $p = V^*MV$, appending $V$ with the necessary monomials.  Let $q$ be such that $q^*Sq=p$, and write $q= \sum_\alpha q_\alpha X^\alpha$.  Let $\hat \alpha$ be maximal, with respect to lexicographical ordering, among all $\alpha$ such that $q_\alpha \neq 0$.

We have $V^*MV=p=q^*Sq=q^*( \sum_{|\alpha|=d} X^{\alpha^*}X^\alpha)q=\sum_{|\alpha|=d} (X^{\alpha}q)^*(X^\alpha q)$. For each $\alpha$, write $X^\alpha q= Q_\alpha V$, where $Q_\alpha$ is the row vector of the coefficients of $X^\alpha q$.  Forming the matrix $Q$ whose rows are the $Q_\alpha$ we get $V^*MV=p=V^*Q^*QV$, hence $V^*(M-Q^*Q)V=0$.  

The polynomials $X^\alpha q$ form a linearly independent set, and in fact have the distinct leading terms $q_{\hat \alpha}X^{\alpha \circ \hat \alpha}$.  It follows that the last $(2n)^{d+deg(q)}$ columns of $Q$ form a block of rank at least $(2n)^d$.  Writing $Q$ in block form $Q=\begin{bmatrix}A&B\end{bmatrix}$ where $B$ is a $(2n)^d\times (2n)^{d+deg(q)}$ matrix, we compute

\[p=V^*Q^TQV=V^*\begin{bmatrix} A^T \\ B^T \end{bmatrix} \begin{bmatrix}A&B\end{bmatrix}V=V^*\begin{bmatrix}A^TA & A^TB\\ B^TA & B^TB\end{bmatrix}V.\]

The $V$ above includes all monomials up to degree $(2n)^{d+deg(q)}$.  Since $V^*(M-Q^*Q)V=0$, we know from the lemma that $M$ cannot differ from $Q^*Q$ in its $(2n)^{d+deg(q)}\times (2n)^{d+deg(q)}$ lower right block; this block equals $B^TB$.  Therefore $M$, an arbitrary matrix representation for $p$, has rank at least $(2n)^d$.\endproof

Alternatively, one might ask whether a Pfister's Theorem holds for products of the usual form. Consider what it would take for $q^*qS$ to be a SOS.  Because $q^*q$ is symmetric, we note that since SOS are symmetric we must have $q^*qS=(q^*qS)^*=Sq^*q$, so that $q^*q$ and $S$ commute.  Since we evaluate these polynomials on tuples of matrices, it is tempting to treat them as symmetric matrices.  In particular, one might guess that if two of them commute, then they are both polynomials in a third polynomial. This happens to be true, and it follows from the following more general theorem from combinatorics:

\begin{theorem}
(Bergman's Centralizer Theorem) Let $K$ be a field, and $K\langle X\rangle$ the ring of polynomials over $K$ in noncommuting variables $X_1,\ldots,X_n$.  Then the centralizer of a nonscalar element in $K\langle X\rangle$ is isomorphic to $K[t]$ for a single variable $t$.
\end{theorem}
The proof is a bit lengthy and can be found in \cite{LO}.  It uses the fact that such a centralizer is integrally closed in its field of fractions together with an easier result in the formal series setting:

\begin{theorem}(Cohn's Centralizer Theorem) Let $K$ be a field and $K\langle \langle X \rangle \rangle$ the ring for formal power series over $K$ in noncommuting variables $X_1, \ldots,X_n$. Then the centralizer of a nonscalar element in $K\langle \langle X \rangle \rangle$ is isomorphic to $K[t]$ for a single variable $t$.

\end{theorem}

These theorems apply despite the superficial difference that we are working with indeterminates $X_1,\ldots,X_n,X_1^*, \ldots,X_n^*$ for which $(X_i^*)^*=X_i$; there are no \emph{polynomial} relations among them, and so we can regard them as $2n$ noncommuting variables $Y_1, \ldots,Y_{2n}$. Armed with Theorem 3.2, we are ready to give the counterexample:

\begin{theorem}

If $p \in \mathbb{R}\langle X, X^*\rangle$, a matrix positive polynomial, is of the form $q^*qS$ with $S=\sum_{|\alpha|=d} X^{\alpha^*}X^\alpha$, then $\text{rank}(p) \geq (2n)^d$. 
\end{theorem}
\proof We will use the previous Theorem 3.1 together with Bergman's Centralizer Theorem.  The main difficulty lies in showing that under the hypotheses, $q^*q$ is actually a polynomial in $S$.

Invoking the centralizer theorem we write $q^*q=f(h(X,X^*))$ and $S=g(h(X,X^*))$ for $h(X,X^*) \in \mathbb{R}\langle X,X^*\rangle$ and $f(t),g(t) \in \mathbb{R}[t]$. It follows from the equation $S=g(h(X,X^*))$ that $g$ must have degree $1$.  To see why, write 

\[h(X,X^*)=c_1X^{\alpha_1}+\ldots+c_lX^{\alpha_l}+(\text{\emph{lower degree terms}}), \hspace{1cm}g(t)=a_kt^k+\ldots+a_0\]

with $c_j,a_i \in \mathbb{R}$. We note that each term $X^{\alpha_{j_1}}\cdots X^{\alpha_{j_k}}$ is symmetric since it must be one of the monomials $X^{\alpha^*}X^\alpha$ in $S$. Supposing $k>1$, we have always that $\alpha_{j_1}=\alpha_{j_k}^*$. This implies that there is just one $\alpha_j$, which is certainly not the case.  Therefore $\text{deg}(g)=1$ and we write $g(t)=at+b$ so that $S=g(h(X,X^*))=ah(X,X^*)+b$ or $h(X,X^*)=1\slash a(S-b)$.

Now we have $q^*q=f(1\slash a (S-b))=r(S)$ for some polynomial $r(t)\in \mathbb{R}[t]$.  Since $r(S)$ has rank equal to 1(it can be expressed as a single noncommutative square), it follows that $r(t)$ is of even degree.  If not, write $r(t)=r_{2k+1}t^{2k+1}+\ldots +r_0$ with $r_{2k+1}\neq0$. Then $r(S)=r_{2k+1}S^{2k+1}+(\text{\emph{lower degree terms}})$ and we have by Theorem 3.1 that $S^{2k+1}=S^kSS^k$ and therefore $r(S)$ itself has rank at least $(2n)^d>1$, a contradiction.  Finally, $tr(t)$ has odd degree and therefore another application of Theorem 3.1 lets us conclude that $p=Sr(S)$ has rank at least $(2n)^d$.\endproof

\bibliographystyle{amsplain}
\bibliography{mynewbib}

\end{document}